\documentclass[reqno,12pt]{article}
\usepackage{amsfonts,latexsym,amscd,amssymb}
\usepackage{amsmath}
\usepackage{epsfig}

\newtheorem{thm}{Theorem}

\begin{document}
\title{\bf Experimental Determination of \\ Ap\'{e}ry-Like Identities for
$\zeta(2n+2)$}

\author{David H.~Bailey\thanks{Lawrence Berkeley Lab, Berkeley, CA 94720 USA.
Email: \texttt{dhbailey@lbl.gov}. Bailey's work is supported by the
Director, Office of Computational and Technology Research, Division
of Mathematical, Information, and Computational Sciences of the U.S.
Department of Energy, under contract number DE-AC02-05CH11231.},
Jonathan M.~Borwein\thanks{Faculty of Computer Science, Dalhousie
University, Halifax, Nova Scotia, B3H 2W5 Canada. Email:
\texttt{jborwein@cs.dal.ca}. Borwein's work is funded by NSERC and
the Canada Research Chair Program.}, and David
M.~Bradley\thanks{Department of Mathematics \& Statistics,
University of Maine, 5752 Neville Hall, Orono Maine, 04469-5752
U.S.A.  Email: \texttt{bradley@math.umaine.edu}}}

\date{\today}

\maketitle

\abstract{We document the discovery of two generating functions for
$\zeta(2n+2)$, analogous to earlier work for $\zeta(2n+1)$ and
$\zeta(4n+3)$, initiated by Koecher and pursued further by Borwein,
Bradley and others.}

\section{Introduction}

Stimulated by recent work in the arena of Ap\'ery-like sums
\cite{apery,sigsam,bbk} we decided to methodically look for series
acceleration formulas for the Riemann zeta function involving
central binomial coefficients in the denominators.   Using the PSLQ
integer relation algorithm, as described below, we uncovered several
new results.  In particular, we document the discovery of two
generating functions for $\zeta(2n+2)$, analogous to earlier work
for $\zeta(2n+1)$ and $\zeta(4n+3)$, initiated by Koecher and
pursued further by Borwein, Bradley and others.  As a conclusion to
a very satisfactory experiment, we have been able to use the
Wilf-Zeilberger technique to prove our results.

An integer relation detection algorithm accepts an $n$-long vector
$\vec x$ of real numbers and a bound $A$ as input, and either
outputs an $n$-long vector $\vec a$ of integers such that the dot
product $a_1 x_1 + \cdots + a_n x_n = 0$ to within the available
numerical precision, or else establishes that no such vector of
integers of length less than $A$ exists.  Here the length is the
Euclidean norm $(a_1^2 + a_2^2 + \cdots + a_n^2)^{1/2}$, derived
from the usual Euclidean metric on $\mathbf{R}^n$.  Helaman
Ferguson's PSLQ~\cite{ferguson, bailey2000} is currently the most
widely used integer relation detection algorithm~\cite[pp.\
230--235]{borwein2004a}, although variants of the so-called LLL
algorithm~\cite{lll} are also commonly employed. Such algorithms
underlie the ``Recognize'' and ``identify'' commands in the
respective computer algebra packages {\sc Mathematica} and {\sc
Maple}. They also play a fundamental role in the investigations we
discuss here.

This origins of this work lay in the existence of infinite series
formulas involving central binomial coefficients in the denominators
for the constants $\zeta(2), \, \zeta(3)$, and $\zeta(4)$.  These
formulas, as well the role of the formula for $\zeta(3)$ in Ap\'ery's
proof of its irrationality, have prompted considerable effort during
the past few decades to extend these results to larger integer
arguments. The formulas in question are
\begin{eqnarray}\label{z2}
\zeta(2) &=& 3 \sum_{k=1}^{\infty} \frac{1}{k^2\binom{2k}{k}}, \\
\smallskip\nonumber\\\label{z3}
\zeta(3) &=& \frac{5}{2} \sum_{k=1}^{\infty} \frac{(-1)^{k+1}}
   {k^3\binom{2k}{k}}, \\
\smallskip \nonumber\\\label{z4}
\zeta(4) &=& \frac{36}{17} \sum_{k=1}^{\infty} \frac{1}
   {k^4 \binom{2k}{k}}.
\end{eqnarray}
Identity (\ref{z2}) has been known since the 19th century---it
relates to $\arcsin^2(x)$---while (\ref{z3}) was variously
discovered in the last century and (\ref{z4}) was noted by
Comtet~\cite[p.\ 89]{comtet}, see~\cite{bbk,vdp}. Indeed, in
\cite{bbk} a coherent proof of all three was provided in the
course of a more general study of such central binomial series and
so-called \emph{multi-Clausen sums}.

These results led many to conjecture that the constant
$\mathcal{Q}_5$ defined by the ratio
\begin{eqnarray*}
\mathcal{Q}_5 &:=& \zeta(5) \bigg/ \sum_{k=1}^{\infty}
   \frac{(-1)^{k+1}}
   {k^5\binom{2k}{k}}
\end{eqnarray*}
is rational, or at least algebraic. However, integer relation
computations using PSLQ and 10,000-digit precision have established
that if $\mathcal{Q}_5$ is a zero of a polynomial of degree at most 25
with integer coefficients, then the Euclidean norm of the vector of
coefficients exceeds $1.24 \times 10^{383}$.  Similar computations for
$\zeta(5)$ have yielded a bound of $1.98 \times 10^{380}$.

These computations lend credence to the belief that $\mathcal{Q}_5$
and $\zeta(5)$ are transcendental.  If algebraic, they almost
certainly satisfy no simple polynomial of low degree. In particular,
if there exist relatively prime integers $p$ and $q$ such that
\begin{eqnarray*}
 \zeta(5) &=& \frac{p}{q} \sum_{k=1}^\infty
  \frac{(-1)^{k+1}}{k^5\binom{2k}{k}},
\end{eqnarray*}
then $p$ and $q$ must be astronomically large. Moreover, a study
of polylogarithmic ladders in the golden ratio
produced~\cite{ag,bbk}
\begin{eqnarray}
   2\zeta(5) -\sum_{k=1}^\infty\frac{(-1)^{k+1}}{k^5 {2k \choose k}}
   &=& \frac 52{\rm Li}_5(\rho)-\frac52{\rm Li_4} (\rho)
    \log \rho + \zeta(3)  \log ^2 \rho -\frac {1}{3}\zeta(2)\log^3
    \rho -\frac1{24}\log^5 \rho, \nonumber \\
&&    \label{rhoeqn}
\end{eqnarray}
where $\rho=(3-\sqrt{5})/2$ and ${\rm Li}_n(z)=\sum_{k=1}^\infty
{z^k}/{k^n}$ is the \emph{polylogarithm} of order $n$. Since the
terms on the right hand side of~(\ref{rhoeqn}) are almost
certainly algebraically independent~\cite{bt}, we see how unlikely
it is that ${\mathcal Q}_5$ is rational. Although the
irrationality of $\zeta(5)$ has not yet been confirmed, it is
known that one of $\zeta(5),\zeta(7),\zeta(9),\zeta(11)$ is
irrational~\cite{Zudilin}.

Given the negative result from PSLQ computations for
$\mathcal{Q}_5$, the authors of \cite{apery} systematically
investigated  the possibility of a multi-term identity of this
general form for $\zeta(2n+1)$. The following were recovered
early~\cite{apery,sigsam} in experimental searches using
computer-based integer relation tools:
\begin{eqnarray}\label{z5}
   \zeta(5) &=& 2\sum_{k=1}^\infty \frac{(-1)^{k+1}}{k^5{2k\choose k}}
    - \frac{5}{2}\sum_{k=1}^\infty \frac{(-1)^{k+1}} {k^3{2k\choose k}}
      \sum_{j=1}^{k-1}\frac{1}{j^2},\\
   \zeta(7) &=& \frac{5}{2}\sum_{k=1}^\infty
               \frac{(-1)^{k+1}}{k^7 {2k \choose k}}
   + \frac{25}{2} \sum_{k=1}^\infty \frac{(-1)^{k+1}}{k^3 {2k \choose k}} \sum_{j=1}^{k-1}\frac{1}{j^4}
     \label{z7} \\
     \zeta(9) &=& \frac{9}{4}\sum_{k=1}^\infty \frac{(-1)^{k+1}}{k^9 {2k
  \choose k}} - \frac{5}{4} \sum_{k=1}^\infty
  \frac{(-1)^{k+1}}{k^7 {2k \choose k}}
  \sum_{j=1}^{k-1}\frac{1}{j^2}
  + 5 \sum_{k=1}^\infty \frac{(-1)^{k+1}}{k^5 {2k
  \choose k} }  \sum_{j=1}^{k-1}\frac{1}{j^4}\nonumber \\\
  && \hspace*{1cm} + \frac{45}{4}\sum_{k=1}^\infty \frac{(-1)^{k+1}}{k^3 {2k
  \choose k} }  \sum_{j=1}^{k-1}\frac{1}{j^6}
  -\frac{25}{4} \sum_{k=1}^\infty \frac{(-1)^{k+1}}{k^3 {2k
  \choose k} }
  \sum_{j=1}^{k-1}\frac{1}{j^4}\sum_{j=1}^{k-1}\frac{1}{j^2},
   \label{z9}\\
  \zeta(11) &=& \frac{5}{2}\sum_{k=1}^\infty \frac{(-1)^{k+1}}{k^{11} {2k
  \choose k}}
  + \frac{25}{2} \sum_{k=1}^\infty
  \frac{(-1)^{k+1}}{k^7 {2k \choose k}} \sum_{j=1}^{k-1}\frac{1}{j^4}
  \nonumber \\
  &-& \frac{75}{4} \sum_{k=1}^\infty
  \frac{(-1)^{k+1}}{k^3 {2k \choose k}} \sum_{j=1}^{k-1}\frac{1}{j^8}
  + \frac{125}{4} \sum_{k=1}^\infty
  \frac{(-1)^{k+1}}{k^3{2k \choose k}} \sum_{j=1}^{k-1}
  \frac{1}{j^4}\sum_{i=1}^{k-1}\frac{1}{i^4}. \label{z11}
\end{eqnarray}
The general formula
\begin{eqnarray}\label{koecher}
   \sum_{k=1}^{\infty }{\frac {1}{k( k^2-x^2)}} &=&
   \frac 12\sum _{k=1}^{\infty }
   \frac{(-1)^{k+1}}{k^3 {2k\choose k}} \frac{5k^2-x^2}{k^2-x^2}
   \prod_{m=1}^{k-1}\left(1-{\frac {{x}^{2}}{{m}^{2}}}\right)
\end{eqnarray}
was obtained by Koecher~\cite{koecher} following techniques of
Knopp and Schur.  It gives~(\ref{z3}) as its first term
and~(\ref{z5}) as its second term but more complicated expressions
for $\zeta(7)$, $\zeta(9)$ and $\zeta(11)$
than~(\ref{z7}),~(\ref{z9}) and~(\ref{z11}).
 The corresponding
result that gives~(\ref{z3}), (\ref{z7}) and~(\ref{z11}) for its
first three terms was worked out by Borwein and
Bradley~\cite{apery}.

Using bootstrapping and an application of the ``Pade'' function
(which in both {\sc Mathematica} and {\sc Maple} produces Pad\'{e}
approximations to a rational function satisfied by a truncated
power series) produced the following remarkable and unanticipated
result~\cite{apery}:
\begin{eqnarray}
\label{z43}
   \sum_{k=1}^\infty \frac{1}{k^3(1-x^4/k^4)} &=&
   \frac{5}{2} \sum_{k=1}^\infty \frac{(-1)^{k+1}}{k^3 {2k \choose
  k} (1-x^4/k^4)} \prod_{m=1}^{k-1} \left(\frac{1+4x^4/m^4}{1-x^4/m^4}
 \right). %= \sum_{n=0}^\infty \zeta(4n+3)x^{4n}.
\end{eqnarray}

The equivalent hypergeometric formulation of~(\ref{z43}) is
\begin{eqnarray*}
   &&
   {}_5F_4\left(\begin{array}{ccccc}2,1+x,1-x,1+ix,1-ix\\
   2+x,2-x,2+ix,2-ix\end{array}\bigg|1\right)\\
   &&\qquad   = \; \bigg(\frac54\bigg)\,{}_6F_5\left(\begin{array}{cccccc}
     2,2, 1+x+ix,1+x-ix, 1-x+ix, 1-x-ix\\
     3/2, 2+x,2-x,2+ix,2-ix\end{array}\bigg|-\frac14\right).
\end{eqnarray*}

The identity~(\ref{z43}) generates~(\ref{z3}), (\ref{z7})
and~(\ref{z11}) above, and more generally gives a formula for
$\zeta(4n+3)$, which for $n>1$ contains fewer summations than the
corresponding formula generated by~(\ref{koecher}).   The task of
proving~(\ref{z43}) was reduced in~\cite{apery} to that of establishing
any one of a number of equivalent finite combinatorial identities.  One
of these latter identities is
\begin{eqnarray}\label{finite3}
   \sum_{k=1}^n{\frac{2n^2}{k^2}}\prod_{i=1}^{n-1}(4k^4+i^4)\bigg/
   \prod_{\substack{i=1 \\ i\neq k}}^{n}(k^4-i^4) &=& {\binom{2n} n}.
\end{eqnarray}
This was proved in~\cite{ag}, so~(\ref{z43}) is an established theorem.
It is now known to be the $x=0$ case of the even more general formula
\begin{eqnarray}\label{bradley}
   \sum_{k=1}^{\infty }{\frac{k}{k^4-x^2k^2-y^4}} &=& \frac12\sum_{k=1}^{\infty }\frac{(
-1)^{k+1}}{{k} {2k \choose k} }\frac{5k^2-x^2}{k^4-x^2k^2-y^4}
\prod_{m=1}^{k-1}{\frac{(m^2-x^2)^2+4y^4}{m^4-x^2m^2-y^4}},
\end{eqnarray}
in which setting $y=0$ recovers~(\ref{koecher}).  The bivariate
generating function identity~(\ref{bradley}) was conjectured by
Henri Cohen and proved by Bradley~\cite{bradley}. It was
subsequently and independently proved by Rivoal~\cite{rivoal}.

Following an analogous---but more deliberate---experimental-based
procedure, as detailed below, we provide a similar general formula
for $\zeta(2n+2)$ that is pleasingly parallel to~(\ref{z43}). It is:

\begin{thm}\label{thm:apery2} Let $x$ be a complex number not
equal to a non-zero integer.  Then
\begin{eqnarray}\label{apery2}
   \sum_{k=1}^\infty \frac{1}{k^2-x^2}
   &=& 3\sum _{k=1}^{\infty } \frac{1}{k^2\,{2k\choose k}
    \left(1- {x}^{2}/k^2 \right)}\prod _{m=1}^{k-1}\left({\frac
{1-4\,{x}^{2}/{m}^{2}}{ 1-{x}^{2}/{m}^{2}}}\right).
\end{eqnarray}
\end{thm}

Note that the left hand side of~(\ref{apery2}) is trivially equal
to
\begin{eqnarray}\label{apery22}
    \sum_{n=0}^{\infty }\zeta(2n+2) x^{2n} &=& \frac{1-\pi x\cot(\pi
 x)}{2x^2}.
\end{eqnarray}
Thus,~(\ref{apery2}) generates an Ap\'ery-like formulae for
$\zeta(2n)$ for every positive integer $n$.

In Section 2 we shall outline the discovery path, and then in
Section 3 we prove~(\ref{apery2})---or rather the equivalent finite
form
\begin{eqnarray}\label{finite2}
  {}_3F_2\left(\begin{array}{ccc} 3n,n+1,-n\\ 2n+1, n+1/2\end{array}
   \bigg|\frac14\right) &=& \frac{ \binom{2n}{n}}{\binom{3n}{n}},
   \qquad 0\le n\in\mathbf{Z}.
\end{eqnarray}

In Section 4 we provide another generating function for which the
leading term or ``seed" is Comtet's formula (\ref{z4}) for
$\zeta(4)$, while the prior generating functions have seeds
(\ref{z2}) and (\ref{z3}).

The paper concludes with some remarks concerning our lack of success
in obtaining formulas analogous to~(\ref{z43}) and~(\ref{apery2})
which would generate the simplest known Ap\'ery-like formulae for
$\zeta(4n+2)$ and $\zeta(4n+1)$, respectively. In this light we
record
\begin{eqnarray}\label{z0}
\sum_{k=1}^{\infty} \frac{1} {\, {2k \choose k}} &=&\frac{2\pi\sqrt{3} +9}{27}, \\
\smallskip\nonumber\\\label{z1}
\sum_{k=1}^{\infty} \frac{(-1)^{k+1}}
   {k\, {2k \choose k}}&=&\frac{2}{\sqrt{5}}\,\log\left(\frac{\sqrt 5+1}{2}
   \right),
\end{eqnarray}
which are perhaps more appropriate seeds in these cases, see
\cite[pp.\ 384--86]{agm}.

\section{Discovering Theorem~\ref{thm:apery2}}

As indicated, we have applied a more disciplined experimental approach
to produce an analogous generating function for $\zeta(2n+2)$. We
describe this process of discovery in some detail here, as the general
technique appears to be quite fruitful and may well yield
 results in other settings.

We first conjectured that $\zeta(2n+2)$ is a rational combination of
terms of the form
\begin{eqnarray}\label{messdown}
   \sigma(2r;[2a_1,\cdots,2a_N]) &:=& \sum_{k=1}^\infty
   \frac{1}{k^{2r} {2k \choose k }} \prod_{i=1}^N \, \sum_{n_i=1}^{k-1}
    \frac1{n_i^{2a_i}},
\end{eqnarray}
  where $r + \sum_{i=1}^N a_i = n+1$, and the
$a_i$ are listed in nonincreasing order (note that the
right-hand-side value is independent of the order of the $a_i$).
This dramatically reduces the size of the search space, while in
addition the sums~(\ref{messdown}) are relatively easy to compute.

One can then write
\begin{eqnarray}\label{guessgf}
   \sum_{n=0}^\infty\zeta(2n+2)\,x^{2n}
   &\stackrel{?}{=}& \sum_{n=0}^\infty
   \sum_{r=1}^{n+1}\sum_{\pi \in
 \Pi(n+1-r)}{\alpha}(\pi)\, \sigma(2r;2\pi)\,x^{2n},
\end{eqnarray}
where $\Pi(m)$ denotes the set of all additive partitions of $m$
if $m>0$, $\Pi(0)$ is the singleton set whose sole element is the
null partition $[\,]$, and the coefficients $\alpha(\pi)$ are
complex numbers. In principle $\alpha(\pi)$ in~(\ref{guessgf})
could depend not only on the partition $\pi$ but also on $n$.
However, since the first few coefficients appeared to be
independent of $n$, we found it convenient to assume that the
generating function could be expressed in the form given above.

For positive integer $k$ and partition $\pi =
(a_,a_2,\ldots,a_N)$ of the positive integer $m$, let
\begin{eqnarray*}
   \widehat{\sigma}_k(\pi) &:=& \prod_{i=1}^N\sum_{n_i=1}^{k-1}
 \frac{1}{n_i^{2a_i}}.
\end{eqnarray*}
Then
\begin{eqnarray*}
   \sigma(2r;2\pi) &=& \sum_{k=1}^\infty
   \frac{ \widehat{\sigma}_k(\pi)}{k^{2r} {2k \choose k }},
\end{eqnarray*}
and from~(\ref{guessgf}), we deduce that
\begin{eqnarray}
  \sum_{n=0}^\infty\zeta(2n+2)\,x^{2n}
  &=& \sum_{n=0}^\infty \sum_{r=1}^{n+1}\sum_{\pi\in\Pi(n+1-r)}
     \alpha(\pi)\sigma(2r;2\pi) x^{2n}\nonumber\\
   &=&\sum_{k=1}^\infty \frac{1}{{{2k}\choose k}} \sum_{r=1}^\infty
\frac{x^{2r-2}}{k^{2r}}\sum_{n=r-1}^{\infty} \sum_{\pi \in
 \Pi(n+1-r)}{\alpha}(\pi)\, \widehat{\sigma}_k(\pi)x^{2(n+1-r)}\nonumber\\
   &=&\sum_{k=1}^\infty \frac{1}{{{2k}\choose k}(k^2-x^2)}
       \sum_{m=0}^{\infty}x^{2m} \sum_{\pi \in \Pi(m)}{\alpha}(\pi)\,
       \widehat{\sigma}_k(\pi) \nonumber\\
  &=& \sum_{k=1}^\infty \frac{1}{{{2k}\choose k}(k^2-x^2)}\,P_k(x)
\end{eqnarray}
where
\begin{eqnarray}
   P_k(x)
  &:=& \sum_{m=0}^{\infty}  x^{2m} \; \sum_{\pi \in
  \Pi(m)}{\alpha}(\pi)\, \widehat{\sigma}_k(\pi), \label{Pk}
\end{eqnarray}
whose closed form is yet to be determined.  Our strategy, as in
the case of~(\ref{z43})~\cite{sigsam}, was to compute $P_k(x)$
explicitly for a few small values of $k$ in a hope that these
would suggest a closed form for general $k$.

Some examples we produced are shown below.  At each step we
``bootstrapped'' by assuming that the first few coefficients of
the current result are the coefficients of the previous result.
Then we found the remaining coefficients (which are in each case
unique) by means of integer relation computations.

In particular, we computed high-precision (200-digit) numerical values
of the assumed terms and the left-hand-side zeta value, and then
applied PSLQ to find the rational coefficients.  In each case we
``hard-wired'' the first few coefficients to agree with the
coefficients of the preceding formula.  Note below that in the sigma
notation, the first few coefficients of each expression are simply the
previous step's terms, where the first argument of $\sigma$
(corresponding to $r$) has been increased by two.

 These initial terms
(with coefficients in bold) are then followed by terms with the other
partitions as arguments, with all terms ordered lexicographically by
partition (shorter partitions are listed before longer partitions, and,
within a partition of a given length, larger entries are listed before
smaller entries in the first position where they differ; the integers
in brackets are nonincreasing):

\begin{eqnarray}
\zeta(2) &=& 3 \sum_{k=1}^\infty \frac {1}{{2 k \choose k} k^2}
  \; = \; 3 \sigma (2,[0]), \\
\zeta(4) &=& \mathbf{3} \sum_{k=1}^\infty \frac{1}{{2 k \choose k} k^4}
  - 9 \sum_{k=1}^\infty \frac{\sum_{j=1}^{k-1} j^{-2}}{{2 k \choose k} k^2}
  \; = \; \mathbf{3} \sigma (4,[0]) - 9 \sigma (2,[2])  \end{eqnarray}
 \begin{eqnarray}
\zeta(6) &=& \mathbf{3} \sum_{k=1}^\infty \frac {1}{{2 k \choose k}
k^6}
 - \mathbf{9} \sum_{k=1}^\infty \frac{\sum_{j=1}^{k-1} j^{-2}}
  {{2 k \choose k} k^4} - \frac{45}{2} \sum_{k=1}^\infty
  \frac{\sum_{j=1}^{k-1} j^{-4}}{{2 k \choose k} k^{2}} \nonumber \\
&& + \frac{27}{2} \sum_{k=1}^\infty
  \sum_{j=1}^{k-1} \frac{\sum_{i=1}^{k-1} i^{-2}}
  {j^2 {2 k \choose k} k^2}, \\
  &=& \mathbf{3} \sigma (6,[]) - \mathbf{9} \sigma (4,[2])
 - \frac{45}{2} \sigma (2,[4]) + \frac{27}{2} \sigma (2,[2,2])\\
\zeta (8) &=& \mathbf{3} \sigma(8,[]) - \mathbf{9} \sigma(6,[2])
 - \mathbf{\frac{45}{2}} \sigma(4,[4])
 + \mathbf{\frac{27}{2}} \sigma (4,[2,2])
 - 63 \sigma (2,[6]) \nonumber \\
 && + \frac{135}{2} \sigma (2,[4,2]) - \frac{27}{2} \sigma (2,[2,2,2]) \\
\zeta (10) &=& \mathbf{3} \sigma (10,[]) - \mathbf{9} \sigma (8,[2])
 - \mathbf{\frac{45}{2}} \sigma (6,[4])
 + \mathbf{\frac{27}{2}} \sigma (6,[2,2])
 - \mathbf{63} \sigma (4,[6]) \nonumber \\
 && + \mathbf{\frac{135}{2}} \sigma (4,[4,2])
 - \mathbf{\frac{27}{2}} \sigma (4,[2,2,2])
 - \frac{765}{4} \sigma (2,[8]) + 189 \sigma (2,[6,2]) \nonumber \\
 && + \frac{675}{8} \sigma (2,[4,4]) - \frac{405}{4} \sigma (2,[4,2,2])
 + \frac{81}{8} \sigma (2,[2,2,2,2]).
\end{eqnarray}

Next from the above results, one can immediately read that
$\alpha([\,]) = 3$, $\alpha([1]) = -9$, $\alpha([2]) = -45/2$,
$\alpha([1,1]) = 27/2$, and so forth.  Table 1 presents the values
of $\alpha$ that we obtained in this manner.

\begin{table}\begin{center}
\begin{tabular}{|l|l|l|l|l|l|}
\hline
Partition & Alpha & Partition & Alpha & Partition & Alpha \\
\hline
 [empty] & 3/1 &  1 & -9/1 &  2 & -45/2 \\
 1,1 & 27/2 &  3 & -63/1 &  2,1 & 135/2 \\
 1,1,1 & -27/2 &  4 & -765/4 &  3,1 & 189/1 \\
 2,2 & 675/8 &  2,1,1 & -405/4 &  1,1,1,1 & 81/8 \\
 5 & -3069/5 &  4,1 & 2295/4 &  3,2 & 945/2 \\
 3,1,1 & -567/2 &  2,2,1 & -2025/8 &  2,1,1,1 & 405/4 \\
 1,1,1,1,1 & -243/40 &  6 & -4095/2 &  5,1 & 9207/5 \\
 4,2 & 11475/8 &  4,1,1 & -6885/8 &  3,3 & 1323/2 \\
 3,2,1 & -2835/2 &  3,1,1,1 & 567/2 &  2,2,2 & -3375/16 \\
 2,2,1,1 & 6075/16 &  2,1,1,1,1 & -1215/16 &  1,1,1,1,1,1 & 243/80 \\
 7 & -49149/7 &  6,1 & 49140/8 &  5,2 & 36828/8 \\
 5,1,1 & -27621/10 &  4,3 & 32130/8 &  4,2,1 & -34425/8 \\
 4,1,1,1 & 6885/8 &  3,3,1 & -15876/8 &  3,2,2 & -14175/8 \\
 3,2,1,1 & 17010/8 &  3,1,1,1,1 & -1701/8 &  2,2,2,1 & 10125/16 \\
 2,2,1,1,1 & -6075/16 &  2,1,1,1,1,1 & 729/16 &  1,1,1,1,1,1,1 & -729/560 \\
 8 & -1376235/56 &  7,1 & 1179576/56 &  6,2 & 859950/56 \\
 6,1,1 & -515970/56 &  5,3 & 902286/70 &  5,2,1 & -773388/56 \\
 5,1,1,1 & 193347/70 &  4,4 & 390150/64 &  4,3,1 & -674730/56 \\
 4,2,2 & -344250/64 &  4,2,1,1 & 413100/64 &  4,1,1,1,1 & -41310/64 \\
 3,3,2 & -277830/56 &  3,3,1,1 & 166698/56 &  3,2,2,1 & 297675/56 \\
 3,2,1,1,1 & -119070/56 &  3,1,1,1,1,1 & 10206/80 &  2,2,2,2 & 50625/128 \\
 2,2,2,1,1 & -60750/64 &  2,2,1,1,1,1 & 18225/64 &  2,1,1,1,1,1,1 & -1458/64 \\
 1,1,1,1,1,1,1,1 & 2187/4480 &&&& \\
\hline
\end{tabular}
\caption{Alpha coefficients found by PSLQ computations}
\end{center}\end{table}

Using these values, we then calculated series approximations to the
functions $P_k(x)$, by using formula (\ref{Pk}).  We obtained:
\begin{eqnarray*}
P_3(x) &\approx& 3 - \frac{45}{4} x^2 - \frac{45}{16} x^4 -
 \frac{45}{64} x^6 - \frac{45}{256} x^8 -
 \frac{45}{1024} x^{10} - \frac{45}{4096}  x^{12}-\frac{45}{16384}
 x^{14} \nonumber \\
 && - \frac{45}{65536} x^{16} \\
P_4(x) &\approx& 3 - \frac{49}{4} x^2 + \frac{119}{144} x^4
 + \frac{3311}{5184} x^4 + \frac{38759}{186624} x^6
 + \frac{384671}{6718464} x^8
\nonumber \\
 &&   + \frac{3605399}{241864704} x^{10}  + \frac{33022031}{8707129344} x^{12}
 + \frac{299492039}{313456656384} x^{14} \\
P_5(x) &\approx& 3 - \frac{205}{16} x^2 + \frac{7115}{2304} x^4
 + \frac{207395}{331776} x^6 + \frac{4160315}{47775744} x^8
 + \frac{74142995}{6879707136} x^{10} \nonumber \\
 && + \frac{1254489515}{990677827584}x^{12}
 + \frac{20685646595}{142657607172096} x^{14}
 + \frac{336494674715}{20542695432781824} x^{16} \\
P_6(x) &\approx& 3 - \frac{5269}{400} x^2 + \frac{6640139}{1440000}
x^4
 + \frac{1635326891}{5184000000} x^6 - \frac{5944880821}{18662400000000}
 x^8 \nonumber \\
 && - \frac{212874252291349}{67184640000000000} x^{10}
 - \frac{141436384956907381}{241864704000000000000} x^{12} \nonumber \\
 && - \frac{70524260274859115989}{870712934400000000000000} x^{14}
 - \frac{31533457168819214655541}{3134566563840000000000000000} x^{16}
\end{eqnarray*}
\begin{eqnarray*}
P_7(x)&\approx& 3 - \frac{5369}{400} x^2 + \frac{8210839}{1440000}
x^4
 - \frac{199644809}{5184000000} x^6
 - \frac{680040118121}{18662400000000} x^8 \nonumber \\
 && - \frac{278500311775049}{67184640000000000} x^{10}
 - \frac{84136715217872681}{241864704000000000000} x^{12} \nonumber \\
 && - \frac{22363377813883431689}{870712934400000000000000} x^{14}
 - \frac{5560090840263911428841}{3134566563840000000000000000} x^{16}.
\end{eqnarray*}

With these approximations in hand, we were then in a position to
attempt to determine closed-form expressions for $P_k(x)$.  This
can be done by using either ``Pade'' function in either {\sc
Mathematica} or {\sc Maple}.  We obtained the following:
\begin{eqnarray*}
P_1(x) &\stackrel{?}{=}& 3 \\ P_2(x) &\stackrel{?}{=}& \frac{3 (4
x^2 - 1)}{(x^2 - 1)}
\\ P_3(x) &\stackrel{?}{=}& \frac{12 (4 x^2 - 1)}{(x^2 - 4)} \\
P_4(x) &\stackrel{?}{=}& \frac{12 (4 x^2 - 1)(4 x^2 - 9)}
         {(x^2 - 4)(x^2 - 9)} \\
P_5(x) &\stackrel{?}{=}& \frac{48 (4 x^2 - 1)(4 x^2 - 9)}
         {(x^2 - 9)(x^2 - 16)} \\
P_6(x) &\stackrel{?}{=}& \frac{48 (4 x^2 - 1)(4 x^2 - 9)(4 x^2 -
25)}
      {(x^2 - 9)(x^2 - 16)(x^2 - 25)} \\
P_7(x) &\stackrel{?}{=}& \frac{192 (4 x^2 - 1)(4 x^2 - 9)(4 x^2 -
25)}
      {(x^2 - 16)(x^2 - 25)(x^2 - 36)}
\end{eqnarray*}

These results immediately suggest that the general form of a generating
function identity is:
\begin{eqnarray}\label{zeta2}
   \sum_{n=0}^\infty\zeta (2 n + 2) {x}^{2n} &\stackrel{?}{=}&
  3\sum_{k=1}^\infty \frac{1}{{2 k \choose k} (k^2-x^2)}
  \prod_{m=1}^{k-1} \frac{4 x^2 - m^2}{x^2 - m^2},
\end{eqnarray}
which is equivalent to~(\ref{apery2}).

We next confirmed this result in several ways:
\begin{enumerate}
\item We  symbolically computed the power series coefficients of
the LHS and the RHS of (\ref{zeta2}), and have verified that they
agree up to the term with $x^{100}$.

\item We verified that $\mathcal{Z}(1/6)$, where $\mathcal{Z}(x)$ is the RHS
of (\ref{zeta2}), agrees with $18 - 3 \sqrt{3} \pi$, computed using
(\ref{apery22}), to over 2,500 digit precision; likewise for
$\mathcal{Z}(1/2)=2, \; \mathcal{Z}(1/3) = 9/2 - 3 \pi/(2 \sqrt{3}), \;
\mathcal{Z}(1/4) = 8 - 2\pi$ and $\mathcal{Z}(1/\sqrt{2}) = 1 - \pi/
\sqrt{2} \cdot \cot(\pi / \sqrt{2})$.

\item We then affirmed that the formula (\ref{zeta2}) gives the
same numerical value as (\ref{apery22}) for the 100 pseudorandom
values $\{m \pi \}$, for $1 \leq m \leq 100$, where $\{\cdot\}$
denotes fractional part.
\end{enumerate}

Thus, we were certain that (\ref{apery2}) was correct and it
remained only to find a proof of Theorem 1.

\section{Proof of Theorem~\ref{thm:apery2}}

By partial fractions,
\begin{eqnarray*}
   \frac{1}{1-z^2/k^2}\prod_{m=1}^{k-1} \frac{1-4z^2/m^2}{1-z^2/m^2}
   &=& \sum_{n=1}^k \frac{c_n(k)}{1-z^2/n^2},
\end{eqnarray*}
where
\begin{eqnarray*}
   c_n(k) &=& \prod_{m=1}^{k-1} (1-4n^2/m^2) \bigg/\prod_{\substack{
   m=1\\ m\ne n}}^k (1-n^2/m^2)
\end{eqnarray*}
if $1\le n\le k$, and $c_n(k) = 0$ if $n>k$ or if $k\ge 2n+1$. It
follows that
\begin{eqnarray*}
   \sum_{k=1}^\infty
   \frac{1}{k^2\binom{2k}{k}(1-z^2/k^2)}\prod_{m=1}^{k-1}\frac{1-4z^2/m^2}{1-z^2/m^2}
   &=& \sum_{k=1}^\infty \frac{1}{k^2\binom{2k}{k}}\sum_{n=1}^k
   \frac{c_n(k)}{1-z^2/n^2}\\
   &=& \sum_{n=1}^\infty \frac{1}{1-z^2/n^2}\sum_{k=n}^{2n}\frac{c_n(k)}{k^2\binom{2k}{k}}.
\end{eqnarray*}
The interchange of summation order is justified by absolute
convergence.  To prove~\eqref{apery2}, it obviously suffices to show
that
\[
   3\sum_{k=n}^{2n} \frac{c_n(k)}{k^2\binom{2k}{k}} =
   \frac{1}{n^2} \quad\Longleftrightarrow\quad S_n:=\sum_{k=n}^{2n} \frac{3}{{2k\choose k}} \prod_{m=1}^{k-1}(4n^2-m^2)
   \bigg/\prod_{\substack{m=1\\m\ne n}}^k (n^2-m^2)=1
\]
for each positive integer $n$. But,
\begin{eqnarray*}
 % &\sum_{k=n}^{2n} \frac{3}{{2k\choose k}} \prod_{m=1}^{k-1}(4n^2-m^2)
 %  \bigg/\prod_{\substack{m=1\\m\ne n}}^k (n^2-m^2)\\
   S_n  &=& \frac{6n^2}{(2n)!}\sum_{k=n}^{2n}\frac{1}{{2k\choose
   k}}\prod_{m=1}^{k-1}(4n^2-m^2)\bigg/\prod_{m=n+1}^k (n^2-m^2)\\
    &=& \frac{(3n)!\, n!}{(2n)! \,(2n)!}\,
    {}_3F_2\bigg(\begin{array}{ccc} 3n, n+1, -n \\ 2n+1, n+1/2
    \end{array} \bigg|    \frac14\bigg).
\end{eqnarray*}
Thus, we have reduced the problem of proving~(\ref{apery2}) to
that of establishing the finite identity
\begin{eqnarray}\label{finite2a}
   T(n) &:=& \frac{(3n)!\, n!}{(2n)! \,(2n)!}\,
    {}_3F_2\bigg(\begin{array}{ccc} 3n,n+1,-n \\ 2n+1, n+1/2 \end{array}
     \bigg|    \frac14\bigg)=1, \qquad n\in\mathbf{Z}^{+}.
\end{eqnarray}
But {\sc Maple} readily simplifies $T(n+1)/T(n)=1$, and since
$T(0)=1$, the identity~(\ref{finite2a}) and hence (\ref{finite2})
and (\ref{apery2}) are established.  If a certificate is desired, we
can employ the Wilf-Zeilberger algorithm.  In {\sc Maple} 9.5 we set
\begin{equation}\label{questionable}
   r:=\frac{{2n\choose n}}{{3n\choose n}},
   \quad
   f:= \frac{(3n+k-1)!\, (n+k)!\, (-n-1+k)!\, (2n)!\, (n-1/2)!\,
  (1/4)^k}{(3n-1)!\, n!\, (-n-1)!\, (2n+k) !\, (n-1/2+k)!\, k!}.
\end{equation}
{\sc Maple} interprets the latter in terms of the Pochammer symbol
\begin{eqnarray*}
   (a)_k &:=& \prod_{j=1}^k (a+j-1)
\end{eqnarray*}
as
\begin{eqnarray*}
   f &=& \frac{(3n)_k \, (n+1)_k \, (-n)_k}{(2n+1)_k \,
   (n+1/2)_k}\cdot \frac{(1/4)^k}{k!},
\end{eqnarray*}
so despite the appearance of~(\ref{questionable}) the issue of
factorials at negative integers does not arise for non-negative
integers $k$ and $n$. Now execute:

\begin{verbatim}
     >  with(SumTools[Hypergeometric]):
     >  WZMethod(f,r,n,k,'certify'): certify;
\end{verbatim}
 which returns the certificate
 \begin{verbatim}
                 /    2                      \
                 \11 n  + 1 + 6 n + k + 5 k n/ k
               - -------------------------------
                  3 (n - k + 1) (2 n + k + 1) n

\end{verbatim}
This proves that summing $f(n,k)$ over $k$ produces $r(n)$, as
asserted.

Indeed,  the (suppressed) output of `WZMethod' is
 the \emph{WZ-pair} $(F,G)$ such that
\begin{eqnarray*}
 F(n+1,k) - F(n,k) &=&
G(n,k+1) - G(n,k),
\end{eqnarray*}
where $F(n,k):= f(n,k)/r(n)$ for $r(n) \neq 0$ and is $f(n,k)$
otherwise. Sum both sides over $k \in \textbf{Z}$ and use the fact
that by construction $G(n,k)
 \to 0$ as $k \to \pm \infty$.
The certificate is
\begin{eqnarray*}
R(n,k) &:=& \frac{G(n,k)}{F(n,k)}.
\end{eqnarray*}

\hfill {\bf QED}

\section{An Identity for $\zeta(2n+4)$}
We compare~(\ref{apery2}) to a result due to
Leshchiner~\cite{lesh} which is stated incorrectly in~\cite{ag},
and which, as the authors say, has a different flavor: for complex
$x$ not an integer,
\begin{eqnarray}\label{lesh}
 \frac 12\sum _{k=1}^{\infty }
  \frac{1}{k^2 {2k\choose k}} \frac{3k^2+x^2}{k^2-x^2}
  \prod_{m=1}^{k-1}\left(1-\frac{x^2}{m^2}\right)
  &=& \sum_{n=1}^{\infty }\frac{(-1)^{n-1}}{n^2-x^2}
  \; = \; \frac{\pi}{2x\sin(\pi x)}- \frac1{2x^2}.
\end{eqnarray}
Using the methods of the previous section---but using a basis of
sums over simplices not hypercubes--- we have likewise now
obtained for complex $x$ not an integer,
\begin{eqnarray}\label{apery4}
   \sum_{k=1}^{\infty }\frac{1}{k^2 {2k\choose k}(k^2-x^2)}
         \prod_{m=1}^{k-1}\left(1-\frac{x^2}{m^2}\right)
   &=& \frac{\pi}{4x^3\sin(\pi x)}-\frac1{x^4}+\frac{3\cos(\pi
   x/3)}{4x^4}.
\end{eqnarray}

To see this, let
\begin{eqnarray*}
   W(x) &:=&
   \bigg(\frac{1/2}{1-x^2}\bigg)\, {}_5F_4\left(\begin{array}{ccccc}
   1,1+x,1+x,1-x,1-x\\
   2,3/2,2+x,2-x\end{array}\bigg|\frac14\right)
\end{eqnarray*}
denote the left hand side of~(\ref{apery4}), and let
\begin{eqnarray*}
   V(x) &:=& \bigg(\frac12\bigg)\, {}_3F_2\left(\begin{array}{ccc} 1, 1+x,
   1-x\\
      2, 3/2 \end{array}\bigg|\frac14\right)
   = \sum_{k=1}^\infty \frac{1}{k^2 {2k\choose
   k}}\prod_{m=1}^{k-1}\bigg(1-\frac{x^2}{m^2}\bigg)\\
   &=&  \sum_{n=1}^\infty (-x^2)^{n-1} \frac{(2
   \arcsin(1/2))^{2n}}{(2n)!}
   = \frac{1-\cos(\pi x/3)}{x^2}.
\end{eqnarray*}
Expanding Leshchiner's series~(\ref{lesh}) now gives
\begin{eqnarray*}
   \frac32\, V(x) + 2x^2\,W(x) &=& \frac{\pi}{2x\sin(\pi x)}-
   \frac1{2x^2}.
\end{eqnarray*}
Solving for $W(x)$ gives~(\ref{apery4}) as claimed.

 \medskip

To recapitulate, we have

\begin{thm} Let $x$ be a complex number, not an integer.
Then
\begin{eqnarray}\label{zeta4}
  \sum_{k=1}^{\infty }\frac{1}{k^2 {2k\choose k}(k^2-x^2)}
         \prod_{m=1}^{k-1}\left(1-\frac{x^2}{m^2}\right)
  &=& {\frac {\pi x\csc( \pi x)+3\cos(\pi x / 3)-4 }{4{x}^{4}}}.
\end{eqnarray}
If  $0\le |x|<1$, then the Maclaurin series for the left hand side
of~(\ref{zeta4}) is equal to
\begin{eqnarray}
 \sum_{n=2}^{\infty }{\frac { \left( -1 \right) ^{n} \left\{
 {3}^{1-2n} -2 B_{2n}
 \left( {2}^{2n-1}-1 \right)  \right\} {\pi }^{2\,n}}{4\,
 \left( 2n \right) !}}{x}^{2n-4} && \nonumber \\
  && \hspace*{-5.5cm} = \; \frac {17}{36}\,\zeta(4)
  + {\frac {313}{648}}\, \zeta(6)\,{x}^{2}+
 {\frac {23147}{46656}}\,\zeta(8)\,{x}^ {4} + {\frac
 {1047709}{2099520}}\,\zeta(10)\,{x}^{6}+\cdots, \nonumber
\end{eqnarray}
where the rational coefficients $B_{2n}$ refer to the even indexed
Bernoulli numbers generated by
\begin{eqnarray*}
   x \coth(x) &=& \sum_{n=0}^\infty \frac{(2x)^{2n}}{(2n)!} B_{2n}.
\end{eqnarray*}
\end{thm}

\medskip

Note that the constant term recaptures (\ref{z4}) as desired---as
taking the limit on the right side  of (\ref{zeta4}) confirms.
Correspondingly, the constant term in  (\ref{lesh}) yields
(\ref{z2}). The coefficient of $x^2$ is
\begin{eqnarray*}
\frac{313}{648}\, \zeta(6) &=& \sum_{k=1}^\infty
\frac{1}{k^6{2k\choose k}}-\, \sum_{k=1}^\infty \frac{1}
{k^4{2k\choose k}} \sum_{j=1}^{k-1}\frac{1}{j^2}  .
\end{eqnarray*}

This all suggests that there should be a unifying formula for our
two identities---(\ref{apery2}) and (\ref{zeta4})---as there is for
the odd cases, see (\ref{bradley}).

\section{Conclusion}

We believe that this general experimental procedure will ultimately
yield results for many other classes of arguments, such as for
$\zeta(4n+m)$, $m=0,1$, but our current experimental results are
negative.

\begin{itemize}
\item[\bf 1.] Considering $\zeta(4n+1)$, for $n=2$ the simplest
evaluation we know is~(\ref{z9}). This is one term shorter than that
given by Rivoal~\cite{rivoal}, which comes from taking the
coefficient of $x^2y^4$ in~(\ref{bradley}).

\item[ \bf 2.]  For $\zeta(2n+4)$ (and $\zeta(4n)$) starting with
(\ref{z4}) which we recall:
\begin{eqnarray*}
\zeta(4) &=& \frac{36\cdot 1}{17} \sum_{k=1}^\infty \frac{1}{
k^{4}\,{2k\choose k}},
\end{eqnarray*}
the  identity for $\zeta(6)$ most susceptible to bootstrapping is
\begin{eqnarray}
\label{z6} \zeta(6) &=& \frac{36\cdot8}{163}\left[\sum_{k=1}^\infty
\frac{1}{k^6{2k\choose k}}+\frac32\, \sum_{k=1}^\infty \frac{1}
{k^2{2k\choose k}} \sum_{j=1}^{k-1}\frac{1}{j^4} \right].
\end{eqnarray}
For $\zeta(8), \zeta(10)$ we have enticingly found:
\begin{eqnarray}
\zeta(8)
 &=&  \frac{36\cdot64}{1373} \left [\sum_{k=1}^\infty \frac{1}
 {k^8 {2k \choose k}}
+ \frac94 \sum_{k=1}^\infty \frac{1}
 {k^4 {2k \choose k}} \sum_{j=1}^{k-1}\frac{1}{j^4}
+ \frac32 \sum_{k=1}^\infty \frac{1}
 {k^2{2k \choose k}} \sum_{j=1}^{k-1}\frac{1}{j^6} \right] \\
\zeta(10) &=&  \frac{36\cdot512}{11143} \left [\sum_{k=1}^\infty
\frac{1}
 {k^{10} {2k \choose k}}
+ \frac{9}{4} \sum_{k=1}^\infty \frac{1} {k^6
 {2k \choose k}} \sum_{j=1}^{k-1}\frac{1}{j^4} + \frac{3}{2} \sum_{k=1}^\infty \frac{1}
 {k^2{2k \choose k}} \sum_{j=1}^{k-1}\frac{1}{j^8}
\right. \nonumber \\
 &&\hspace{1in}  +  \frac{9}{4} \sum_{k=1}^\infty \frac{1}
 {k^4 {2k \choose k}} \sum_{j=1}^{k-1}\frac{1}{j^6} \left.+\frac{27}{8} \sum_{k=1}^\infty \frac{1}
 {k^2 {2k \choose k}} \sum_{j=1}^{k-1}\frac{1}{j^4}\sum_{i=1}^{j-1}\frac{1}{i^4}
 \right].\qquad
\end{eqnarray}
But this pattern is not fruitful; the pattern stops at $n=10$.
\end{itemize}

\end{document}